\documentclass[12pt]{amsart}

\title{Weierstrass models}

\author{Ying Zong}

\address{}
\email{zongying.upenn@gmail.com}

\date{} 

\begin{document}

\maketitle

\section{Introduction}

\smallskip

This paper draws inspiration from lagrangian fibrations on irreducible symplectic manifolds in general, and more specifically, from the question whether all these fibrations have smooth base. Recall that lagrangian fibrations on irreducible symplectic manifolds are equidimensional, at least in the framework of projective varieties, see Matsushita \cite{matsushita}. And when such a base $S$ is also assumed to be smooth, $S$ is then a projective space by Hwang \cite{Hwang}.

\smallskip

In this direction, we proved the following theorem 1 in \cite{equi_regularity}. The theorem 2 is a further effort in studying a certain class of these fibrations.

\smallskip

\smallskip

{\bf Theorem 1.} --- \emph{Let $S$ be a locally noetherian normal algebraic space of residue characteristic zero with regular formal fibers. Assume that there is a morphism $f: X\to S$, which is locally of finite type, equidimensional and admits \'{e}tale local sections around every point of $S$, where $X$ is a regular algebraic space. Then $S$ is regular and $f$ is flat.}

\smallskip

\smallskip

{\bf Corollary.} --- \emph{If a lagrangian fibration $f: X\to S$ on an irreducible symplectic manifold $X$ has local holomorphic sections around a point $s$ of $S$, and if $X$ and $S$ are both projective varieties, then $S$ is smooth over $\mathbf{C}$ at $s$.}

\smallskip

For then $f$ admits formal sections around $s$, and one can apply Theorem 1.

\smallskip

\smallskip

{\bf Definition.} --- \emph{Keeping the hypotheses of Theorem 1, we say that $f$ is a Weierstrass model on $X$ if furthermore }:
\smallskip

\emph{a) $f$ is locally for the \'{e}tale topology of $S$ projective.}

\smallskip

\emph{b) For every geometric codimension $\leq 1$ point $s$ of $S$, the fiber $f^{-1}(s)$ is irreducible.}

\smallskip

\emph{c) Each connected component $S_0$ of $S$ has a geometric point $s_0$ such that $f^{-1}(s_0)$ admits a $k(s_0)$-abelian variety structure.}

\smallskip

The definition is motivated by elliptic $K3$, whose study often benefits from their contraction to plane models, the Weierstrass models. It may happen that such a contraction is trivial, a phenomenon that is generalized here. The following theorem says that Weierstrass models are, to a large extent, group theoretical, being compactification of torsors in a very tight form. We propose a possible application below.

\smallskip

Note that in the condition \emph{c}) if $S_{(s_0)}$ denotes the strict localization of $S$ at $s_0$, then $f\times_SS_{(s_0)}$ is smooth, and $X\times_SS_{(s_0)}$ admits an abelian scheme structure over $S_{(s_0)}$ (\emph{cf.} \cite{basic} 7.2), and hence in particular, for every geometric generic point $t$ of $S_{(s_0)}$, $f^{-1}(t)$ has a $k(t)$-abelian variety structure. 
\smallskip

It follows that, if $f$ is a Weierstrass model on $X$, then $f\times_SS'$ is a Weierstrass model on $X\times_SS'$ for every regular morphism $S'\to S$. 

\smallskip

\smallskip

{\bf Theorem 2.} --- \emph{Let $f: X\to S$ be a Weierstrass model on $X$. Then there exist }:

\smallskip

\emph{a) An $S$-smooth group algebraic space $A$ with connected fibers, which is N\'{e}ron at every codimension $1$ point of $S$.}

\smallskip

\emph{b) An $A$-torsor $Q$ on $S$ for the \'{e}tale topology.}

\smallskip

\emph{c) An open $S$-immersion $j: Q\to X$.}

\smallskip

\emph{Such $(A, Q, j)$ is unique up to unique isomorphisms and of formation compatible with every base change by a regular morphism $S'\to S$.}

\smallskip

\smallskip

We call $Q$ the \emph{identity component}, and $A$ the \emph{albanese}, of the Weierstrass model. As in the case of plane models, $A$ in general has fibers with non-trivial unipotent part. 
\smallskip

The proof, unlike that of Theorem 1, does not need assume $S$ has excellent local rings, and is insensitive to characteristics.
\smallskip

\smallskip

{\bf Proposition 3.} --- \emph{Let $f: X\to S$ be a Weierstrass model on $X$, with $S$ of finite type over a field $k$. Let $s$ be a point of $S$. Assume that, in a neighborhood of $f^{-1}(s)$, $f$ is lagrangian with respect to a symplectic structure $\Psi$ on $X/k$.}
\smallskip

\emph{Then $A_s$ is of abelian rank at least $\mathrm{tr. deg.}\ k(s)/k$.}

\smallskip

\smallskip

This estimate says, for example, that, for $s$ of codimension $1$ in $S$, $A_s$, if it degenerates, is extension of an abelian variety $B$ by a form of $\mathbf{G}_m$ or $\mathbf{G}_a$. In the $\mathbf{G}_m$ case, if $\overline{s}$ denotes a geometric point over $s$, there is a $\mathbf{P}^1$-bundle $P$ over $B_{\overline{s}}$ with two sections $0, \infty: B_{\overline{s}}\to P$ such that $f^{-1}(\overline{s})$ is obtained from $P$ by pinching $0$ and $\infty$; this $f^{-1}(\overline{s})$ is one particular instance of a \emph{minimal model of $L^*$} we classified in \cite{minimal}.

\smallskip

Keeping the hypotheses of Theorem 2, assume now that $X$ is an irreducible symplectic manifold with its symplectic structure $\Psi$ and hyperkahler metric $\omega_g$ (\emph{cf}. \cite{beauville}). Then $H^0(\mathcal{O}_S)=H^0(f_{*}\mathcal{O}_X)=H^0(\mathcal{O}_X)=\mathbf{C}$, $S$ is proper smooth over $\mathbf{C}$, and $f$ is a lagrangian fibration, at least when $S$ is projective over $\mathbf{C}$ (\cite{matsushita} Theorem 1). Notice that all terms of the exact sequence
\[0\to f^*\Omega^1_{S}|Q\to \Omega^1_{Q}\to \Omega^1_{Q/S}\to 0\] have natural action by $A$. So by descent one has the exact sequence
\[0\to \Omega^1_{S}\to \omega^1_{Q}\to \omega^1_{Q/S}\to 0\] with $\omega^1_{Q}$ and $\omega^1_{Q/S}$ consisting of the $A$-invariant forms of $\Omega^1_{Q}$ and $\Omega^1_{Q/S}$, respectively. Also, $\Psi|Q$ by descent induces a symplectic form $\psi$ on $\omega^1_{Q}$, relative to which $\Omega^1_{S}$ is lagrangian. 
\smallskip

We now define a hermitian metric on $S$: for each local holomorphic $1$-form $\alpha$ on $S$, let it be considered as an $A$-invariant form on $Q$ and then extended uniquely to $X$ ($\mathrm{codim}(X-Q, X)\geq 2$), say $\alpha'$; let $|\alpha'|_g$ denote its norm relative to the hyperkahler metric $\omega_g$ on $X$. Clearly, the fiberwise integral
\[\alpha\mapsto |\alpha|^2:=\int_{X/S} |\alpha'|^2_g\ \omega^n_g \] where $n$ is the relative dimension of $X/S$, induces a hermitian metric on $S$, which I speculate is Kahler-Einstein given the symmetry it enjoys.

\smallskip

\section{Proof of Theorem 2}

\smallskip

Recall that N\'{e}ron models are of formation compatible with all base change by regular morphisms of spectra of discrete valuation rings \cite{neron} 7.2/1. The question is thus local for the \'{e}tale topology of $S$. We may assume $f$ projective, and that $S$ is a scheme, quasi-compact, irreducible and with generic point $t$.

\smallskip

I) \emph{Case where $S=t$ }:

\smallskip

By assumptions of a Weierstrass model, $X\times_t\overline{t}$ admits a $\overline{t}$-abelian variety structure, for $\overline{t}$ a geometric point of $t=S$. The assertion then follows by \cite{basic} 7.2: $X=Q$ admits a unique torsor structure under its albanese $A$, dual abelian variety of $\mathrm{Pic}^o_{X/t}$.  

\smallskip

II) \emph{Case where $S$ is the spectrum of a discrete valuation ring }:

\smallskip

Let $s$ be the closed point of $S$, and let $U$ consist of all those points of $X$ where $f$ is smooth. Since we assume that $f$ admits \'{e}tale local sections around $s$, $U_s$ is non-empty and hence dense in the geometrically irreducible closed fiber $X_s$. 
\smallskip

We have shown that the generic fiber $U_t=X_t$ is a torsor under its albanese $A_t$. It follows that $U$ is weak $S$-N\'{e}ron \cite{neron} 3.5/1 by definition, and, as $U_s$ is irreducible, is also \emph{minimal}, \emph{loc.cit.} 4.3/2, and finally, by \emph{loc.cit.} 6.5 Corollary 3, that $U$ is $S$-dense open in the $S$-N\'{e}ron model $Q$ of $U_t=X_t$, and that $Q$ is a torsor under the N\'{e}ron model $A$ of $A_t$. 
\smallskip

The fiber $Q_s$, thus $A_s$ as well, is geometrically irreducible, since $U_s$ is geometrically irreducible and dense in $Q_s$. 

\smallskip

To finish, it suffices to show that $U=Q$. 
\smallskip

Let $L$ be an ample invertible module on $X$; such $L$ exists by the hypothesis that $X$ be projective over $S$. Notice that $Q-U$ (resp. $X-U$) is of codimension $\geq 2$ in $Q$ (resp. $X$). As $Q$ is regular, $Q$ is parafactorial along $Q-U$, so $L|U$ uniquely extends to an invertible module $N$ on $Q$. It remains only to see that $N$ is ample on $Q$, because 
\[X=\mathrm{Proj} \sum H^0(L^{\otimes n})=\mathrm{Proj} \sum H^0(L^{\otimes n}|U)=\mathrm{Proj} \sum H^0(N^{\otimes n}).\] That $N$ is ample follows immediately from Raynaud's theory of cube \cite{raynaud} VI 2.1; it is enough to find an integer $n>0$ and a section $\sigma\in H^0(N^{\otimes n})$ such that
\[W:=\{y\in Q, \sigma(y)\neq 0\}\] is contained in $U$, \emph{affine},  with $W_s\neq \emptyset$. Such a section $\sigma$ exists, since 
\[H^0(N^{\otimes n})=H^0(L^{\otimes n}|U)=H^0(L^{\otimes n})\] and $L$ is ample.
\smallskip

III) \emph{General case }:

\smallskip

III.1) \emph{Uniqueness and existence of $A$ }:

\smallskip

\emph{i) If $A$ exists, it is unique up to unique isomorphisms }:

\smallskip

This follows by \cite{deligne} Note (1) 3: The functor $A\mapsto A_t$ is fully faithful for those $A$ that is N\'{e}ron at all codimension $1$ points of $S$.

\smallskip

\emph{ii) In order that the $S$-smooth group algebraic space $A$ exists, it suffices that it exists over an $S$-algebraic space $S'$ that is smooth surjective over $S$ }:

\smallskip

For, by \emph{i)} and I), the asserted $S'$-smooth group algebraic space $A'$ will admit a descent datum relative to $S'\to S$. Such a descent is effective by Artin \cite{artin} 6.3.

\smallskip

\emph{iii) Let $S'$ be the open sub-scheme of $X$ consisting of all points where $f$ is smooth. Then $S'\to S$ is surjective }:

\smallskip

Because $S$ is regular and we assume that $f$ admits \'{e}tale local sections around every point of $S$, \emph{cf.} \cite{neron} 3.1/3.

\smallskip

\emph{iv) For the existence of $A$, one can assume that $f$ admits an $S$-section, say $e: S\to X$ }:

\smallskip

One takes the base change $S'\to S$, with $S'$ as in \emph{iii)}, and applies \emph{ii)}.

\smallskip

\emph{v) Let $e: S\to X$ be as in iv), let $U$ consist of all points of $X$ where $f$ is smooth, and let $V$ consist of all those points $x$ of $U$ such that $x$ and $e(f(x))$ lie in the same connected component of $U_{f(x)}$. Then $V$ is an open sub-scheme of $U$ and has geometrically irreducible fibers over $S$ }:

\smallskip

That $V$ is open in $U$ follows from EGA IV 15.6.5. Its $S$-fibers are by definition connected, thus geometrically connected, since the $S$-section $e: S\to X$ factors through $V$, and finally geometrically irreducible, since $V\to S$ is smooth.

\smallskip

\emph{vi) Let $V$ be as in v). There exists on $V$ an $S$-birational group law with zero section $e:S\to V$. Such a structure is unique up to unique isomorphisms }:

\smallskip

Notice that $V_t=X_t$, which is, in a unique way, an abelian variety $A_t$ with identity $e(t)$. And, we have shown above that, when localized at any codimension $1$ point $s$ of $S$, $V=U$ coincides with the N\'{e}ron model of $A_t$ over $\mathrm{Spec}(\mathcal{O}_{S, s})$. The sought-after $S$-birational composition law $m$ on $V$ therefore, by the technique of ``passage \`{a} la limite projective'', exists over an open sub-scheme $S'$ of $S$, with $S-S'$ of codimension $\geq 2$ in $S$. 
\smallskip

Now \cite{raynaud} IX 1.1 implies that the domain of definition of $m$ is an open neighborhood of $(e\times_Se)(S)$, thus $m$ exists throughout; it verifies the axioms of commutativity, associativity, etc, because it does over $t$. 
\smallskip

\emph{viii) There exists an open neighborhood $V'$ of $e(S)$ in $V$ such that $(V', m|V', e)$ admits an $S$-dense open immersion into an $S$-group algebraic space $A$, compatible with the group laws.}
\smallskip

One applies SGA 3 XVIII 3.2+3.5+3.7 to the birational group law of \emph{vii)}.

\smallskip

\emph{ix) The $S$-group algebraic space $A$ of viii) is a scheme, separated, smooth, of finite type, and with connected fibers over $S$ }:

\smallskip

That $A$ is a scheme (resp. $S$-separated, resp. of finite type over $S$) follows by \cite{neron} 6.2 b), (resp. a), resp. b)), because $V'$ is quasi-projective over $S$. And, $A$ is $S$-smooth with connected fibers, because $V'$ is $S$-smooth with geometrically connected fibers and the morphism $V'\times_SV'\to A$, induced by the group law of $A$, is faithfully flat, which one can verify fiber by fiber.
\smallskip

\emph{x) One has $V=A$ }:

\smallskip

Let $L$ be an $S$-ample invertible module on $X$, which exists since $f$ is by assumption projective. The restriction $L|V'$ extends uniquely to an invertible module $N$ on $A$. Indeed, as $A$, being smooth over $S$, is regular and $\mathrm{codim}(A-V', A)\geq 2$, one has that $A$ is parafactorial along $A-V'$. Observe also that $X-V'$ is of codimension $\geq 2$ in $X$. Now, $N$ is $S$-ample on $A$, which one proves by localizing at each point $s$ of $S$ and argues similarly as in II). This concludes the proof, because
\[X=\mathrm{Proj} \sum f_*L^{\otimes n}=\mathrm{Proj} \sum (f|V')_*(L^{\otimes n}|V')=\mathrm{Proj} \sum a_*N^{\otimes n},\] where $a: A\to S$ denotes the structural morphism of $A$.

\smallskip

III.2) \emph{When $f$ admits an $S$-section $e$, the $(A, Q, j)$ exists as desired and remains so after every base change by a regular morphism }:

\smallskip

One takes $j: Q\to X$ to be the inclusion of $A=V$ in $X$.

\smallskip

III.3) \emph{The $(A, Q, j)$ in III.$2)$ is independent of the choice of $S$-sections of $f$, when such sections exist }:

\smallskip 

Let $e_1, e_2$ be two $S$-sections of $f$, and let $(A_i, Q_i, j_i)$ denote the triple constructed relative to $e_i$, $i=1, 2$.  Identify $Q_i$ with its image in $X$, and write $Q=Q_1\cap Q_2$. 
\smallskip

Observe that $Q_i-Q$ is of codimension $\geq 2$ in $Q_i$. Indeed, $Q_1$ and $Q_2$ coincide over every codimension $\leq 1$ point of $S$, by I)+II). It follows that, this intersection $Q$, when viewed as a birational correspondence between $Q_1$ and $Q_2$, is, according to Weil \cite{neron} 4.4/1, in fact biregular. So the claim $Q_1=Q_2$ is justified. 

\smallskip

Write $+$ for the group law of $A=Q$ with origin $e_1$. There exists a unique element $\alpha_{12}\in A(S)$ satisfying 
\[\alpha_{12}+e_1=e_2.\]

III.4) \emph{Completion of the proof }:

\smallskip

Clearly, varying the sections of $f$ over all $S$-smooth schemes $S'$, we have produced, through the system of $\alpha_{12}$, a unique descent datum gluing the locally constructed $(A, Q, j)$. This finishes the proof.

\smallskip

\section{Proof of Proposition 3}

\smallskip

The field $k$ necessarily is of characteristic $0$, and $X$ and $S$ are smooth over $\mathrm{Spec}(k)$. 
\smallskip

Restricting to a neighborhood of $f^{-1}(s)$, we write $L: T_{X/k}\simeq \Omega^1_{X/k}$ for the isomorphism induced by the symplectic structure $\Psi$. 
\smallskip

Let the closed image of $Q_s$ in $X_s$ be $W$, and let $Z\to W$ be a resolution of singularities, centered in $W-Q_s$, so that $Q_s$ admits naturally an open dense immersion $\varphi$ into $Z$. 
\smallskip

Consider the composition
\[T_{Z/k}\to T_{X/k}\otimes\mathcal{O}_Z\stackrel{L\otimes 1}{\rightarrow}\Omega^1_{X/k}\otimes\mathcal{O}_Z\to\Omega^1_{Z/k}\to \Omega^1_{Z/s}\] where 
\[T_{Z/k}\to T_{X/k}\otimes\mathcal{O}_Z\]
\[\Omega^1_{X/k}\otimes\mathcal{O}_Z\to \Omega^1_{Z/k}\to\Omega^1_{Z/s}\] are induced by the morphism $Z/s\to X/k$.
\smallskip

Call this composition $T_{Z/k}\to \Omega^1_{Z/s}$, $\beta$; we show below that it is of kernel $T_{Z/s}$. Note first that all terms of the exact sequence
\[0\to T_{Z/s}\to T_{Z/k}\to T_{s/k}\otimes\mathcal{O}_Z\to 0\] are locally free over $\mathcal{O}_Z$.
\smallskip

\emph{a) One has $\beta|T_{Z/s}=0$, hence $\beta$ induces a morphism $\overline{\beta}: T_{s/k}\otimes\mathcal{O}_Z\to \Omega^1_{Z/s}$ }:

For, $\beta|T_{Z/s}$, a section over $Z$ of $\mathcal{H}\mathrm{om}(T_{Z/s}, \Omega^1_{Z/s})\simeq (\Omega^1_{Z/s})^{\otimes 2}$, vanishes on $Q_s$ by assumption that $f$ is lagrangian.

\smallskip

\emph{b) The morphism $\overline{\beta}$ is injective }:
\smallskip

It is injective over $Q_s$, because $f$ is lagrangian, and hence is injective throughout, since $T_{s/k}\otimes\mathcal{O}_Z$ is locally free over $\mathcal{O}_Z$.

\smallskip

\emph{c) One has $\mathrm{dim}_{k(s)}\ H^0(Z, \Omega^1_{Z/s})\geq \mathrm{dim}_{k(s)}T_{s/k}=\mathrm{tr. deg.}\ k(s)/k$ }:

\smallskip

This follows by \emph{b)}.
\smallskip

\emph{d) Reduction to a lemma }:
\smallskip

The restriction $\varphi^*: H^0(Z, \Omega^1_{Z/s})\to H^0(Q_s, \Omega^1_{Q_s/s})$ is injective. To finish, one combines \emph{c)} with the lemma below applied to the open immersion $\varphi_{\overline{s}}: Q_{\overline{s}}\hookrightarrow Z_{\overline{s}}$, where $\overline{s}$ is a geometric point over $s$.

\smallskip

\smallskip

{\bf Lemma.} --- \emph{Let $k$ be a field of characteristic $0$. Let $A$ be a connected $k$-algebraic group, with generic point $\eta$, and with Chevalley decomposition \[1\to L\to A\to B\to 1\] so that $B$ is an abelian variety and $L$ affine connected. Let $\varphi$ be a $k$-rational map of $A$ to a proper $k$-algebraic space $Z$.}
\smallskip
 
\emph{Then
\[\mathrm{Image}\ \varphi^*: H^0(Z, \Omega^1_{Z/k})\to \Omega^1_{\eta/k}\] is contained in
\[\mathrm{Image}\ H^0(B, \Omega^1_{B/k})\to H^0(A, \Omega^1_{A/k})\to \Omega^1_{\eta/k}\] In particular, the abelian rank of $A$ verifies :
\[\mathrm{dim}(B)\geq \mathrm{rank}_k(\varphi^*)\]}

\begin{proof} Recall that every algebraic group over a field of characteristic $0$ is smooth (Cartier). In particular, $A$ and $L$ are smooth over $k$, the projection $p: A\to B$ is smooth and faithfully flat, $p^*: H^0(B, \Omega^1_{B/k})\to \Omega^1_{\eta/k}$ is injective and of rank $\mathrm{dim}(B)$, and the lemma amounts to asserting the existence of a unique $k$-linear map 
\[\rho: H^0(Z, \Omega^1_{Z/k})\to H^0(B, \Omega^1_{B/k})\] such that 
\[\varphi^*=p^*\circ \rho.\]

Under this form, we may and will assume $k=\mathbf{C}$.
\smallskip

Write $U$ for the domain of definition of $\varphi$, and consider $\varphi$ as a morphism of $U$ to $Z$. One has $\mathrm{codim}(A-U, A)\geq 2$, since $A$, being $k$-smooth, is normal, and $Z/k$ is proper.

\smallskip

\emph{a) One may assume $\varphi$ is an open schematically dense immersion }:

\smallskip

Decompose $\varphi$ as
\[U\stackrel{\varphi'}{\hookrightarrow}Z'\stackrel{f}{\rightarrow}Z\] with $f$ proper, and $\varphi'$ an open schematically dense immersion. Clearly, the lemma holds for $(\varphi, Z)$ if it holds for $(\varphi', Z')$. It suffices to consider the latter.

\smallskip

\emph{b) One may assume furthermore $Z$ smooth over $k$ }:

\smallskip

Assume by \emph{a)} that $\varphi$ is an open schematically dense immersion. Let $f: Z'\to Z$ be a resolution of singularities, centered in $Z-\varphi(U)$, so that $\varphi$ factors as
\[U\stackrel{\varphi'}{\hookrightarrow}Z'\stackrel{f}{\rightarrow}Z\] with $\varphi'$ an open schematically dense immersion. As in \emph{a)}, it is sufficient to study $(\varphi', Z')$.

\smallskip

\emph{c) The compostion
\[H^0(Z, \Omega^1_{Z/k})\stackrel{\varphi^*}{\rightarrow}H^0(U, \Omega^1_{A/k})\to H^0(U, \Omega^1_{A/B})\] is zero }:

\smallskip

This is obvious if $p: A\to B$ is an isomorphism, i.e., if $L=1$. Assume $L$ non-trivial. 
\smallskip

By structure of affine algebraic groups, it is enough to show that, for each $b\in B(k)$ and each curve $C$ in $p^{-1}(b)\simeq L$ of the form $\mathbf{G}_a$ or $\mathbf{G}_m$, the composition
\[(\varphi|C)^*: H^0(Z, \Omega^1_{Z/k})\stackrel{\varphi^*}{\rightarrow}H^0(U, \Omega^1_{A/k})\to H^0(U, \Omega^1_{A/B})\to H^0(C\cap U, \Omega^1_{C/k})\] is zero. 
\smallskip

Now $Z/k$ being proper, if $C\cap U\neq \emptyset$,  
$\varphi|C: C\cap U\to Z$ extends to a morphism of $\mathbf{P}^1$ to $Z$, but $H^0(\mathbf{P}^1, \Omega^1_{\mathbf{P}^1/k})=0$, hence \emph{c)}.

\smallskip

\emph{d) Let $(\varphi, Z)$ be as in a), b). For each $\alpha\in H^0(Z, \Omega^1_{Z/k})$, one has
\[\varphi^*\alpha=\sum a_i\ p^*e^i\] where $\{e^i\}$ is a fixed basis of $H^0(B, \Omega^1_{B/k})$, and $a_i\in H^0(U, \mathcal{O}_A)$ }:

\smallskip

This follows from \emph{c)} and the exact sequence
\[0\to p^*\Omega^1_{B/k}\to \Omega^1_{A/k}\to \Omega^1_{A/B}\to 0\]

\emph{e) All coefficients $a_i\in H^0(U, \mathcal{O}_A)$ in d) belong to $k$ }:

\smallskip

As $H^0(U, \mathcal{O}_A)=H^0(A, \mathcal{O}_A)$, the claim \emph{e)} is clear if $B=1$ or if $p: A\to B$ is an isomorphism. Assume $n=\mathrm{dim}(B)>0$ and $L$ non-trivial. Choose a volume form
\[\Phi\in H^0(B, \Omega^n_{B/k})\] and let $e_j\in H^0(B, \Omega^{n-1}_{B/k})$ be such that 
\[e^i\wedge e_j=\delta^i_j\ \Phi\] $\forall\ i, j\in \{1,\cdots, n\}$.
\smallskip

One deduces from the expression of $\varphi^*\alpha$ in \emph{d)} that
\[\varphi^*\alpha\wedge p^*e_j=a_j\ p^*\Phi\] and by differentiation that
\[d_{U/k}(a_j)\wedge p^*\Phi=0\] Indeed,
\[d_{Z/k}(\alpha)=0,\ d_{B/k}(\Phi)=0,\ d_{B/k}(e_j)=0\] by the $E_1$-degeneration of Hodge-De Rham spectral sequences. 
\smallskip

Now, $d_{U/k}(a_j)\wedge p^*\Phi=0$ is equivalent to $d_{U/B}(a_j)=0$. In fact, in terms of the Koszul filtration 
\[K^i:=p^*\Omega^i_{B/k}\wedge \Omega^{*-i}_{A/k}\] on $\Omega^*_{A/k}$, one has $d_{U/k}(a_j)\wedge p^*\Phi\in H^0(U, K^n\Omega^{n+1}_{A/k})$, which, under the isomorphism 
\[K^n=K^n/K^{n+1}\simeq \Omega^{*-n}_{A/B}\otimes \Omega^n_{B/k}\] is mapped to $d_{U/B}(a_j)\otimes \Phi$.

\smallskip

And, $d_{U/B}(a_j)=0$ implies that, when viewed as an analytic function,
\[a_j\in \mathrm{Image}\ H^0(V^{an}, \mathcal{O}^{an}_B)\to H^0(U^{an}, \mathcal{O}^{an}_A)\] where $V=p(U)$, since $p: U\to V$ is smooth surjective with geometrically connected fibers. 
\smallskip

Observe that $U$ intersects with all codimension $1$ fibers of $p$. One has hence $\mathrm{codim}(B-V, B)\geq 2$, $H^0(V^{an}, \mathcal{O}^{an}_B)=H^0(B^{an}, \mathcal{O}^{an}_B)=k$. This concludes that $a_j\in k$.
\smallskip  

\emph{f) Completion of the proof }:
\smallskip

By \emph{e)}, all $a_i\in k$. The map 
\[\rho: \alpha \mapsto \sum a_i e^i\] is as desired.

\end{proof}

\smallskip


\bibliographystyle{amsplain}


\end{document}